\documentclass[smallextended,numbook,runningheads]{svjour3}     

\smartqed  

\usepackage{graphicx}

\usepackage{subfigure}
\usepackage{amsmath}
\usepackage{hyperref}

\begin{document}

\title{Solving Tensor Structured Problems with Computational Tensor Algebra}

\titlerunning{Computational Tensor Algebra}

\author{Oleksii V. Morozov \and Patrick R. Hunziker}


\institute{
	Oleksii V. Morozov \at University Hospital of Basel, Petersgraben 4, 4031 Basel, Switzerland \\
	\email{MorozovA@uhbs.ch} 
	\and Patrick R. Hunziker \at University Hospital of Basel, Petersgraben 4, 4031 Basel, Switzerland \\             
    \email{HunzikerP@uhbs.ch}
}

\maketitle

\begin{abstract}
\label{abstract}

Since its introduction by Gauss, \textbf{\textit{Matrix Algebra}} has facilitated understanding of scientific problems, hiding distracting details and finding more elegant and efficient ways of computational solving. Today's largest problems, which often originate from multidimensional data, might profit from even higher levels of abstraction. We developed a framework for \textbf{\textit{solving tensor structured problems with tensor algebra}} that unifies concepts from tensor analysis, multilinear algebra and multidimensional signal processing. In contrast to the conventional matrix approach, it allows the formulation of multidimensional problems, in a multidimensional way, preserving structure and data coherence; and the implementation of automated optimizations of solving algorithms, based on the commutativity of all tensor operations. Its ability to handle large scientific tasks is showcased by a real-world, 4D medical imaging problem, with more than 30 million unknown parameters solved on a current, inexpensive hardware. This significantly surpassed the best published matrix-based approach.


\end{abstract}

\section{Introduction}
\label{Intro}
	Computational problems with tensor structure, which involve large amounts of multidimensional data, arise in many fields of science and engineering \cite{ArigovindanImageRec2005,EldenTensorLS2005,Pereyra1973,Baumgartner2005,Lynch1964,KoldaTensorReview2009}. The standard way of dealing with such problems is not intrinsically multidimensional, and assumes reduction of problem formulation to the classical matrix formalism, which is built on two-dimensional (matrix) and one-dimensional (vector) data structures. This is achieved by the reordering of multidimensional terms from the original problem formulation, into matrices and vectors. The derived matrix formulation is then treated using well-known techniques and algorithms of matrix algebra. Finally the solution is reconstituted into its natural multidimensional form.
	
While this approach reduces a multidimensional problem to the well-known standard form of matrix algebra, it introduces certain limitations. The resulting formulation is no longer explicitly multidimensional, - a fact that may create difficulty in identifying and understanding important properties inherent to the particular problem \cite{Muti2003}. Vectorization of multidimensional objects may lead to loss of spatial data coherence \cite{ComonTensorDecomp2000,WangDataComprCANDECOMP2004}, which can adversely affect the performance of solving algorithms. In many cases, the derived problem formulation does not have a straightforward, intuitive connection with the process of generating efficient solving algorithms.

Currently, a growing interest in objects with more then two dimensions -- tensors, has become evident in the scientific community. Introduced in Tensor Analysis and Multilinear Algebra, tensors gained the attention of practitioners of diverse fields (see the excellent review by Kolda \cite{KoldaTensorReview2009}). However, translating tensor mathematics to a convenient computational framework raises many issues \cite{Martin2004}; including, convenient notation, ease of algorithm implementation, performance.
	
Here we propose a comprehensive framework for \textbf{\textit{solving tensor structured problems by tensor algebra}} that allows natural and elegant formulation of multidimensional problems using multidimensional data structure directly. The proposed framework is based on a generalization of the concepts of matrix algebra to multiple dimensions; it incorporates and unifies existing approaches from multidimensional signal processing \cite{Blahut1985,Guoan2004,AldroubiWaveletsInMedicine} , recent developments in the field of tensor analysis \cite{Heinbockel2001}, and multilinear algebra \cite{Tucker1966,Leibovici1998,LathauwerHOSVD2000,Harshman1970,KoldaTensorReview2009,AcarTensorSurvey2009,Martin2004}. It is based on ordered vector spaces and a generalized definition of tensor multiplication, based on the concept of tensor equations and the tensor inverse, and induces the design of fast computational solving algorithms with built-in expression analysis for automated generation of efficient implementations on serial and parallel computers. In our companion paper we show how to computationally solve large real-world problems using our framework, and claim that this novel approach offers a natural, higher level abstraction to solving a broad range of very large computational problems which arise in multidimensional signal processing and other scientific fields.	

\section{Tensor nomenclature, basic tensor operations}
\label{sec1}
In our framework, multidimensional data and multidimensional transformations are described based on a formalism originating in Tensor Analysis \cite{Heinbockel2001} and Physics \cite{EinsteinFoundRelTheory}. Each object in the formulation represents a multidimensional entity -- a tensor. Tensors describing multidimensional data belong to a space that is the outer product of a number of vector spaces (see \ref{Tensor_notation}). A tensor can be represented by a multidimensional array of its components\footnotemark[1]. \footnotetext[1]{Further in the document the term tensor is identified with tensor components} For example, a scalar field defined on a discrete uniform three-dimensional grid with extents $[N_X, N_Y, N_Z]$ can be expressed as a third order tensor ${\cal{T}}^{x y z}$, where tensor order is determined by the number of indices. Note that in our tensor notation, direct correspondence between designation of tensor indices $(x,y,z)$ and dimensions of the physical space $(X,Y,Z)$ is typically used. Tensors indices can be either covariant (subscripts) or contravariant (superscripts), depending on the way  the tensor components transform in respect to a change of basis (see \ref{cov_contrav}). This property of a tensor index is called variance. For the common case of Euclidean spaces with orthonormal bases, there is no difference between covariant and contravariant indices; thus the choice of index variance is influenced by the context and consistency of tensor expressions.	

Consider a transformation applied to tensor ${\cal{T}}^{x y z}$ which results in a tensor in the same tensor product space. The transformation can be expressed by ${\cal{A}}^{x_1 y_1 z_1}_{x y z}$. Its application to $\cal{T}$ implies the use of classical tensor operations, such as the outer product, which expands tensor order

\begin{equation}
{\cal{A}}^{x_1 y_1 z_1}_{x y z} \cdot {\cal{T}}^{x y z} = {\cal{C}}^{x_1 y_1 z_1 x y z}_{x y z}
\end{equation}

and the contraction, which reduces tensor order

\begin{equation}
{\cal{C}}^{x_1 y_1 z_1 x y z}_{x y z} \equiv \sum_x{\sum_y{\sum_z{{\cal{C}}^{x_1 y_1 z_1 x y z}_{x y z}}}} = 	{\cal{D}}^{x_1 y_1 z_1}
\end{equation}

For designation of the contraction, we use the Einstein convention \cite{EinsteinFoundRelTheory} which assumes implicit summation over a pair of indices with equal designation but different variance. From the practical point of view, it is convenient to combine the outer product and the contraction into a single operation, in analogy with the matrix multiplication. Thus by allowing simultaneous contraction over multiple corresponding index pairs we get

\begin{equation}
{\cal{A}}^{x_1 y_1 z_1}_{x y z} \cdot {\cal{T}}^{x y z} = {\cal{D}}^{x_1 y_1 z_1}
\label{first_product}
\end{equation}

	Our framework relies on the fact that a tensor is an object which can implicitly contain all its properties as indices and component values. Due to this object nature in a practical implementation (\ref{first_product}) can be reduced to D := A*T, where "*" denotes tensor multiplication.

\subsection{Concept of ordered spaces, commutativity of tensor multiplication}
In some cases, multidimensional transformations can be decomposed to a product of one-dimensional transformations. Such transformations are called separable. They are reduced to sequential application of one-dimensional transformations from the decomposition. Using tensor formalism in three dimensions this can be expressed as

\begin{equation}
{\cal{A}}^{x_1 y_1 z_1}_{x y z} \cdot {\cal{T}}^{x y z} = {\cal{B}}^{x_1}_{x} \cdot {\cal{C}}^{y_1}_{y} \cdot {\cal{D}}^{z_1}_{z} \cdot {\cal{T}}^{x y z} = {\cal{E}}^{x_1 y_1 z_1}
\label{separable_transf}
\end{equation}

In conventional tensor formalism (and likewise in matrix formalism), the result of expression (\ref{separable_transf}) depends on the order of multiplicands (lack of commutativity). This is due to concatenation of non contracted indices, which is used for assurance of a unique result of the tensor product \footnotemark[2]. \footnotetext[2]{With operations on mixed tensors, it can be rather difficult to track the overall order of array indices in the results. Thus in some tensor literature a special notation is used, where unique index order is defined by use of "dot" symbol ($E^{.xy}_{z..}$). But this brings an undesirable side-effect in significantly decreasing the readability of tensor expressions.} In contrast, separable data handling operations, like filtering or resampling in signal processing, lead to the same result independently of the order in which separable transformations are applied. To avoid this formal, but practically important mismatch, and to add more flexibility to handling of tensor expressions, we introduce the following constraint: in analogy to physical space, which is ordered (right hand rule), we require vector spaces used for construction of the tensor product to have a unique predefined order. This constraint \footnotemark[3] leads to a fundamental new property of the framework, compared to conventional approaches for matrix and tensors. \textbf{\textit{The tensor multiplication becomes commutative}}, in addition to being associative. Formally this means that indices of the resulted tensor have unique positions determined by a predefined order of vector spaces. Thus complex tensor products can be evaluated in arbitrary order but lead to the same result. For example, for ordered spaces $X,Y,Z$ we can write

\footnotetext[3]{Note that this constraint does not limit the expressiveness of tensors, as the order of vector spaces can be changed if necessary}

\begin{equation}
{\cal{A}}^{x_1}_{x} \cdot {\cal{B}}^{y_1}_{y} \cdot {\cal{C}}^{z_1}_{z} \cdot {\cal{T}}^{x y z} = {\cal{B}}^{y_1}_{y} \cdot {\cal{A}}^{x_1}_{x} \cdot {\cal{C}}^{z_1}_{z} \cdot {\cal{T}}^{x y z} = {\cal{B}}^{y_1}_{y} \cdot  {\cal{C}}^{z_1}_{z} \cdot {\cal{A}}^{x_1}_{x} \cdot {\cal{T}}^{x y z} = {\cal{E}}^{x_1 y_1 z_1}
\end{equation}	

The commutativity of the tensor multiplication, achieved hereby, considerably increases the ease of handling tensor expressions. One of the important practical values of this property consists in simplification of semantic analysis of tensor expressions performed for the purpose of reduction of computational and storage requirements. For example, expression ${\cal{A}}^{z_1 t}_{z} \cdot {\cal{B}}^{y_1}_{y} \cdot {\cal{C}}^{x_1}_{x} \cdot {\cal{T}}^{x y z}_{t}, (N_X < N_Y < N_Z < N_T, N_X N_Y > N_T)$, which according to conventional tensor and matrix formalisms can be evaluated only in a single chain order, in our framework can be computed by reordering ${\cal{C}}^{x_1}_{x} \cdot \left( {\cal{B}}^{y_1}_{y} \cdot \left( {\cal{A}}^{z_1 t}_{z} \cdot {\cal{T}}^{x y z}_{t} \right) \right)$. This gives the best computational performance with minimal memory requirements for storing intermediate results. At this point, it is possible to use the tensor multiplication, in combination with tensor addition and subtraction, defined from multilinearity of tensor product space (see \ref{tensor_operations}) for formulating multidimensional problems in a natural and elegant multidimensional representation. Note that formulations expressed in tensor terms are also convenient for differentiation in respect to multidimensional terms (see \ref{tensor_diff}), an important prerequisite for multidimensional optimization problems.

\subsection{Extension of the tensor notation}
According to Einstein notation it is not permitted to have more than one index with same designation and variance. However, in our settings we do allow this, with agreement on some consistency rules. As an example, consider the following expression

\begin{equation}
{\cal{A}}^{i}_{x} \cdot {\cal{B}}^{i}_{y} = {\cal{H}}^{i}_{xy}
\end{equation}

\emph{This operation is a tensor analogue of the Khatri-Rao product \cite{SmildeMultiwayAnalysis2004} which is a matching elementwise (over $i$) Kronecker product of two sets of vectors.} Note that two formally  equivalent indices are merged into a single one, without dependency on outer product or contraction applied to other indices. But the constraint is, that indices can be contracted only once after all possible merges resulting in a single pair of superscript and subscript. Here is an example which represents the most general form of the tensor product, combining outer product, contraction and elementwise product

\begin{equation}
{\cal{A}}^{i}_{x} \cdot {\cal{D}}^{x}_{i} \cdot {\cal{B}}^{i}_{y} \cdot {\cal{P}}_{iz} \cdot {\cal{W}}_i = \left( {\cal{A}}^{i}_{x} \cdot {\cal{B}}^{i}_{y} \right) \cdot \left( {\cal{D}}^{x}_{i} \cdot {\cal{P}}_{iz} \cdot {\cal{W}}_i \right) = {\cal{H}}^{i}_{xy} \cdot {\cal{N}}^{x}_{iz} = {\cal{R}}_{yz}
\end{equation}

Notice how neatly the proposed notation unifies different types of products which exist separately in Matrix Algebra (Matrix product, Kronecker product, Khatri-Rao product).
 
\section{Tensor equations and tensor inverse}
\label{sec2}
Expression ${\cal{A}}^{x_1 y_1 z_1}_{x y z} \cdot {\cal{U}}^{x y z} = {\cal{B}}^{x_1 y_1 z_1}$ represents a relation between a multidimensional input ${\cal{U}}^{x y z}$ and a multidimensional output ${\cal{B}}^{x_1 y_1 z_1}$. This relation is determined by a transformation ${\cal{A}}^{x_1 y_1 z_1}_{x y z}$. When the input is unknown, the expression describes an inverse problem. In cases where a unique solution exists, the existence of a tensor inverse is implied -- a transformation which maps ${\cal{A}}$ to the identity tensor

\begin{equation}
\tilde{{\cal{A}}}^{x_2 y_2 z_2}_{x_1 y_1 z_1} \cdot  {\cal{A}}^{x_1 y_1 z_1}_{x y z} \cdot {\cal{U}}^{x y z} = {\delta}^{x_2 y_2 z_2}_{x y z} \cdot {\cal{U}}^{x y z} = {\cal{U}}^{x_2 y_2 z_2} = \tilde{{\cal{A}}}^{x_2 y_2 z_2}_{x_1 y_1 z_1} \cdot {\cal{B}}^{x_1 y_1 z_1}
\label{TensorEq}
\end{equation}

where ${\delta}^{x_2 y_2 z_2}_{x y z} = {\delta}^{x_2}_{x} \cdot {\delta}^{y_2}_{y} \cdot {\delta}^{z_2}_{z}$ is the identity transformation in three dimensions expressed in terms of Kronecker delta (see \ref{kronecker_delta}). \emph{Note that equation \ref{TensorEq} can be represented in an equivalent form as a matrix mapping from a vector space to a vector space, with proper reshaping of tensor terms. In cases when system tensor $\cal{A}$ has some structure, treatment of the problem using the tensor representation can be more advantageous than using the conventional matrix formalism, that is explained in the next sections.}

\section{Tensor solvers}
\label{sec4}

\subsection{Direct tensor solvers, iterative tensor solvers, Krylov subspace tensor solvers}

Here we introduce tensor solvers which can be explicitly applied to solve such tensor equations computationally. The most straightforward way to find the solution tensor and, if necessary, a tensor inverse - is based on use of a tensor extension of Gauss elimination or a tensor analogue of LU decomposition - both belonging to the class of direct tensor solvers. These methods can benefit from the explicit tensor structure of equation coefficients, which is typically manifested by handling  multidimensional "sub-blocks" of nonzero coefficients. In many practical problems, however, tensor equations are so large that solving them using direct methods becomes prohibitive. In this instance, we propose the class of iterative tensor solvers. For example, one can use a tensor Jacobi solver which is based on the iterative computation of tensor multiplication and belongs to the class of stationary iterative solvers. The pseudocode of the algorithm and an actual software implementation are shown on Fig.~\ref{Tensor_Jacobi}.

\begin{figure}[]
\begin{center}

\begin{align*}
&{\cal{R}}^{x_1 y_1 z_1} = {\cal{A}}^{x_1 y_1 z_1}_{xyz} \cdot {\cal{U}}^{xyz} - {\cal{B}}^{x_1 y_1 z_1} \\
&\mbox{for $i$=0 until convergence} \\
&\quad {\cal{U}}^{xyz} = {\cal{U}}^{xyz} - {\cal{R}}^{x_1 y_1 z_1} \cdot {\cal{E}}^{xyz}_{x_1 y_1 z_1} \\
&\quad {\cal{R}}^{x_1 y_1 z_1} = {\cal{A}}^{x_1 y_1 z_1}_{xyz} \cdot {\cal{U}}^{xyz} - {\cal{B}}^{x_1 y_1 z_1} \\
&\mbox{end}
\end{align*}
\subfigure[]{}

\begin{align*}\footnotesize
&\mbox{\textbf{PROCEDURE} TensorJacobi(A,B: Tensor;}\\
&\quad \quad \quad \quad \quad \qquad \quad \qquad \qquad \qquad \mbox{threshold: Scalar}\\
&\quad \quad \quad \quad \quad \qquad \quad \qquad \qquad \qquad \mbox{): Tensor;}\\
&\mbox{\textbf{VAR} R, E, U: Tensor;}\\
&\mbox{\textbf{BEGIN}}\\
&\quad \quad \mbox{\bigskip E := InvMainDiag(A,B);}\\
&\quad \quad \mbox{R := A*U - B;}\\
&\quad \quad \mbox{\textbf{WHILE} R+*R $>$ threshold \textbf{DO}}\\
&\quad \quad \quad \mbox{U := U - R*E;}\\
&\quad \quad \quad \mbox{R := A*U - B;}\\
&\quad \quad \mbox{\textbf{END};}\\
&\quad \quad \mbox{\textbf{RETURN} U;}\\
&\mbox{\textbf{END} TensorJacobi;} \\
\end{align*}
\subfigure[]{}

\begin{align*}
&\mbox{\textbf{VAR}} \\
&\quad \mbox{w: World;} \\
&\quad \mbox{A, B, C: Tensor;} \\
&\mbox{\textbf{BEGIN}} \\
&\qquad \qquad \qquad \qquad \qquad \quad \quad \mbox{(* define a world with vector spaces *)} \\
&\quad \mbox{\textbf{NEW}(w);} \\
&\quad \mbox{w.DefineSpace('X',128);} \\
&\quad \mbox{w.DefineSpace('Y',129);} \\
&\quad \mbox{w.DefineSpace('Z',130);} \\
&\qquad \qquad \qquad \qquad \qquad \quad \quad \mbox{(* define system tensor *)} \\
&\quad \mbox{A := Tensors.Laplacian(w,"x\^{}1,x\_,y\^{}1,y\_,z\^{}1,z\_");} \\
&\qquad \qquad \qquad \qquad \qquad \quad \quad \mbox{(* define right hand side *)} \\
&\quad \mbox{NEW(B,w,"x\^{}1,y\^{}1,z\^{}1");} \\
&\quad \mbox{B.Fload("rhs.dat");} \\
&\qquad \qquad \qquad \qquad \qquad \quad \quad \mbox{(* solve the problem *)} \\
&\quad \mbox{C := TensorJacobi(A,B,1.0E-4);} \\
&\quad \cdots \\
&\mbox{\textbf{END}} \\
\end{align*}
\subfigure[]{}

\end{center}
\caption{ (a) Pseudocode of tensor Jacobi iterator for three dimensions (${\cal{E}}$ is the inverse of diagonal tensor constructed from main diagonal of ${\cal{A}}$ ). (b) Actual implementation of solver for any number of dimensions using our tensor library designed in Tensor Oberon programming language, which offers algebraic operators and optimized implementations for the basic tensor operations. Here "*" and "+*" stand for tensor multiplication and inner product respectively. (c) A simple code example of using the solver.}
\label{Tensor_Jacobi} 
\end{figure}
	
	Note that the presented object based implementation of the Jacobi solver \textbf{\textit{does not use explicit indexing}} of tensor components. All information about indices is contained in tensor objects which are defined once at the construction stage. This is sufficient for performing further computations with automatic semantic analysis for optimizing the performance, and without complicating the solving process by explicit use of indices.
	
From the point of view of computational and storage requirements, the presented iterative solver is particularly attractive for application to tensor structured problems. As a basic example, consider the Poisson equation in three dimensions ${\nabla}^{2}{u(x,y,z)} = f(x,y,z)$ which after discretization on a uniform grid can be expressed by the following tensor equation

\begin{equation}
{\cal{A}}^{x_1}_{x} \cdot {\cal{U}}^{x y_1 z_1} + {\cal{B}}^{y_1}_{y} \cdot {\cal{U}}^{x_1 y z_1} + {\cal{C}}^{z_1}_{z} \cdot {\cal{U}}^{x_1 y_1 z} = {\cal{F}}^{x_1 y_1 z_1}
\end{equation}

Here ${\cal{U}}$ is an unknown tensor, ${\cal{F}}$ is the observation tensor, ${\cal{A, B, C}}$ are one-dimensional finite difference transformations which are equivalent to filters with impulse response dependent on the order of the Laplacian approximation. This decomposition allows very efficient computation of the tensor multiplication by performing sequential separable transformations along individual dimensions of the tensor ${\cal{U}}$. The separability renders the computation of the whole set of equation coefficients unnecessary, thus dramatically reducing storage requirements of the presented iterative algorithm (obviously, for such simple example, exploiting the tensor characteristics of the problem is also possible using the standard matrix approach with careful index bookkeeping).		

	Another type of iterative algorithms which are frequently used for solving large inverse problems, is the family of Krylov subspace solvers, which includes GMRES, CG, BICGSTAB and other solvers \cite{SaadIteratSolv2000}. In our extension of this solver class to tensors, most of them can be expressed as an iterative application of the tensor multiplication (see Fig.~\ref{Tensor_CG}) and thus profit from the same benefits discussed above for the tensor Jacobi iterator. In contrast to Jacobi solver which only converges well for a limited class of problems, Krylov solvers perform significantly better for the large majority of problems.
	
	Fig.~\ref{Tensor_CG} presents a pseudocode and Oberon implementation of the tensor Conjugate Gradients (CG) algorithm applied to solve a 3D problem ${\cal{A}}^{x_1 y_1 z_1}_{xyz} \cdot {\cal{C}}^{xyz} = {\cal{B}}^{x_1 y_1 z_1}$.

\begin{figure}[ht]

\begin{center}

\begin{tabular}{l l}
${\cal{R}}^{x_1 y_1 z_1} = {\cal{B}}^{x_1 y_1 z_1} - {\cal{A}}^{x_1 y_1 z_1}_{xyz} \cdot {\cal{C}}^{xyz} \mbox{ \% initial residual}$ & R := B-A*C; \\

${\cal{P}}^{x_1 y_1 z_1} = {\cal{R}}^{x_1 y_1 z_1}$ \% initial search directions & P := R;  \\

 & D := Delta([C.inds, -B.inds]); \\

$\rho = \langle {\cal{R}}^{x_1 y_1 z_1} , {\cal{R}}^{x_1 y_1 z_1} \rangle$ & rho := R+*R; \\

for $i$=1 by 1 until convergence & \textbf{REPEAT} \\

\quad ${\cal{Q}}^{x_1 y_1 z_1} = {\cal{A}}^{x_1 y_1 z_1}_{xyz} \cdot ({\cal{P}}^{x_1 y_1 z_1} \cdot \delta^{xyz}_{x_1 y_1 z_1})$ & \quad Q := A*(P*D); \\

\quad $\alpha = {\rho \over \langle  {\cal{P}}^{x_1 y_1 z_1}, {\cal{Q}}^{x_1 y_1 z_1} \rangle}$ & \quad alpha := rho/(P+*Q); \\

\quad \% update solution approximation & \\

\quad ${\cal{C}}^{xyz} = {\cal{C}}^{xyz} + \alpha \cdot {\cal{P}}^{x_1 y_1 z_1} \cdot \delta^{xyz}_{x_1 y_1 z_1}$ & \quad C := C + alpha*P*D; \\

\quad ${\cal{R}}^{x_1 y_1 z_1} = {\cal{R}}^{x_1 y_1 z_1} - \alpha \cdot {\cal{Q}}^{x_1 y_1 z_1}$ \% update residual & \quad R := R - alpha*Q; \\

\quad ${\rho}_1 = \langle {\cal{R}}^{x_1 y_1 z_1} , {\cal{R}}^{x_1 y_1 z_1} \rangle$ & \quad rho1 = R+*R; \\

\quad \% update search directions & \\

\quad ${\cal{P}}^{x_1 y_1 z_1} = {\cal{R}}^{x_1 y_1 z_1} + { \left( {\rho}_1 \over \rho \right)} \cdot {\cal{P}}^{x_1 y_1 z_1}$ & \quad P := R + (rho1/rho)*P; \\

\quad $\rho = {\rho}_1$ & \quad rho := rho1; \\

end & \textbf{UNTIL} rho $>$ threshold \\
\end{tabular}

\end{center}
\caption{Pseudocode and given side by side actual implementation of Conjugate Gradient iterator}
\label{Tensor_CG}
\end{figure}

\subsection{Tensor Multigrid algorithms unify algebraic multigrid and geometric multiresolution} In the matrix domain, multigrid methods for solving linear system of equations are known for their fast convergence and computational efficiency. These methods perform iterative improvement of a solution approximation, based on smoothing of the solution error on multiple scales of the problem. Basically, there are two types of multigrid/multiresolution strategies: \textbf{\textit{geometric}} multiresolution preserves spatial coherence in multidimensions by working on datasets at different resolutions but its formulation is often problem specific. \textbf{\textit{Algebraic}} matrix multigrid algorithms (AMG) offer a 'black box' approach for solving matrix equations, but they rely on values of matrix coefficients only and risk neglecting data coherence, at a potential cost of degraded performance.

	We have used tensor equations to express both approaches in a unified fashion. The proposed tensor multigrid algorithm (TMG) has an important advantage over AMG: in contrast to the matrix algorithm, which neglects multidimensionality and the structure of the problem, TMG unifies both concepts of geometric and algebraic multiresolution, by preserving and exploiting multidimensionality and spatial data coherence (see Fig.~\ref{TMG}).
	
\begin{figure}[ht]
\includegraphics[width=1.0\textwidth]{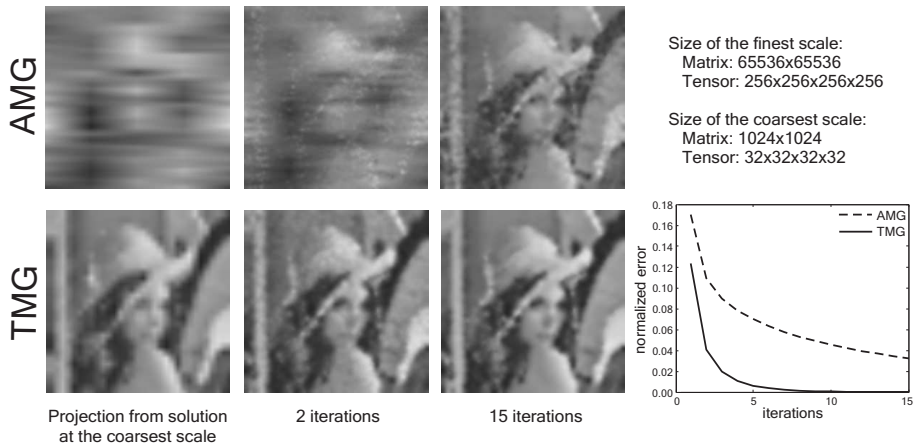}
\caption{Performance of matrix (AMG) and tensor (TMG) multigrid algorithms in application to a problem of image reconstruction from an incomplete set of points. The tensor approach (TMG), through preservation of data locality in the coarser scales, yields higher quality solutions and faster convergence.}
\label{TMG_vs_AMG} 
\end{figure}

	Fig.~\ref{TMG} shows pseudocode for a Tensor Algebraic Multigrid V-cycle algorithm applied to solve a 3D tensor problem ${\cal{A}}^{x_1 y_1 z_1}_{xyz} \cdot {\cal{C}}^{xyz} = {\cal{B}}^{x_1 y_1 z_1}$. Problems in higher dimensions are handled in a similar way. Fig.~\ref{TMG_vs_AMG} compares visually how AMG and TMG are applied to a problem of reconstruction of a two-dimensional image from a few arbitrarily taken samples, using non-uniform spline interpolation \cite{ArigovindanImageRec2005}. This example shows that performance of the AMG algorithm is significantly degraded due to data distortion, which is introduced at coarse matrix scales, whilst the TMG converges to a good approximation of the solution within a few iterations, due to its preservation of spatial coherence.

\begin{figure}[ht]

\begin{center}

\begin{align*}
& \mbox{function $C^{xyz}$ = TMGVCycle($scale$,${\cal{A}}^{x_1 y_1 z_1}_{xyz}$,${\cal{C}}^{xyz}$,${\cal{B}}^{x_1 y_1 z_1}$)} \\
& \quad \mbox{if $scale \neq last$ then}\\
& \quad \quad {\cal{C}}^{xyz} = \mbox{Presmooth(${\cal{A}}^{x_1 y_1 z_1}_{xyz}$,${\cal{C}}^{xyz}$,${\cal{B}}^{x_1 y_1 z_1}$)}\\
&\quad \qquad \qquad \qquad \qquad \qquad \qquad \qquad \quad \quad \mbox{\% \textit{get preliminary constructed system scale and}} \\
&\quad \qquad \qquad \qquad \qquad \qquad \qquad \qquad \quad \quad \mbox{\% \textit{corresponding restriction/prolongation operators}} \\
& \quad \quad \left[ \hat{{\cal{A}}}^{u_1 v_1 w_1}_{uvw}, \hat{{\cal{C}}}^{uvw}, \hat{{\cal{B}}}^{u_1 v_1 w_1}, {\cal{P}}^{u_1}_{x_1}, \hat{{\cal{P}}}^{x}_{u}, {\cal{S}}^{v_1}_{y_1}, \hat{{\cal{S}}}^{y}_{v}, {\cal{T}}^{w_1}_{z_1}, \hat{{\cal{T}}}^{z}_{w} \right] = \mbox{GetScale($scale$) }\\
&\quad \qquad \qquad \qquad \qquad \qquad \qquad \qquad \quad \quad \mbox{\% \textit{separable restriction of the residual}} \\
& \quad \quad \hat{{\cal{B}}}^{u_1 v_1 w_1} = {\cal{P}}^{u_1}_{x_1} \cdot {\cal{S}}^{v_1}_{y_1} \cdot {\cal{T}}^{w_1}_{z_1} \cdot \left( {\cal{B}}^{x_1 y_1 z_1} - {\cal{A}}^{x_1 y_1 z_1}_{xyz} \cdot {\cal{C}}^{xyz} \right)\\
& \quad \quad \hat{{\cal{C}}}^{uvw} = 0 \\
& \quad \quad \hat{{\cal{C}}}^{uvw} = \mbox{TMGVCycle($scale+1$,$\hat{{\cal{A}}}^{u_1 v_1 w_1}_{uvw}$,$\hat{{\cal{C}}}^{uvw}$,$\hat{{\cal{B}}}^{u_1 v_1 w_1}$) \% \textit{V-cycle on next scale}}\\
& \quad \quad {\cal{C}}^{xyz} = {\cal{C}}^{xyz} + \hat{{\cal{P}}}^{x}_{u} \cdot \hat{{\cal{S}}}^{y}_{v} \cdot \hat{{\cal{T}}}^{z}_{w} \cdot \hat{{\cal{C}}}^{uvw} \mbox{ \% \textit{coarse grid correction}}\\
& \quad \quad {\cal{C}}^{xyz} = \mbox{Postsmooth(${\cal{A}}^{x_1 y_1 z_1}_{xyz}$,${\cal{C}}^{xyz}$,${\cal{B}}^{x_1 y_1 z_1}$)}\\
& \quad \mbox{else} \\
& \quad \quad  {\cal{C}}^{xyz} = \mbox{DirectSolve(${\cal{A}}^{x_1 y_1 z_1}_{xyz}$,${\cal{B}}^{x_1 y_1 z_1}$) }\\
& \quad \mbox{end} \\
& \mbox{end}
\end{align*}

\end{center}
\caption{Pseudocode of a TMG V-cycle algorithm for three dimensions}
\label{TMG}
\end{figure}	 

\section{Automatic tensor expression analysis for high-performance implementations} 
\label{sec3}
Because the high level tensor formulation contains the complete tensorial structure information of a complex mathematical problem, automatic expression analysis can be performed, resulting in automated software optimization and code generation. We have experimentally verified this aspect in our tensor library and have found that such automatic analysis is capable of optimizing code, with regard to computational and memory requirements, in both a problem-specific and hardware-specific manner. This result profits strongly from the commutativity of all tensor operations – the significant differentiating feature over standard matrix formulations which are non-commutative with respect to matrix multiplication. Some aspects of such automated code generation have also been discussed in \cite{Baumgartner2005}. This particular topic is a promising area for future research.

\section{Application of the framework to large real-world problems}
\label{sec4} 
Our framework was successfully applied to solve a large real-world computational problem –- a spline-based global variational reconstruction \cite{ArigovindanImageRec2005} of a multidimensional signal from incomplete and spatially scattered measurements in four dimensions. The problem was originally formulated in terms of matrix algebra. In four dimensions, the matrices involved in the computations become extremely large, essentially prohibiting computation of the problem on current standard hardware. For example, for a moderate size of the reconstruction grid $128 \times 128 \times 128 \times 16$ we have to deal with a matrix which is represented by 7'514'144'528 nonzero elements (in compressed block-band diagonal storage, ~30 GBytes in single precision). Using our framework, we reformulated the problem in terms of tensors. The resulting tensor formulation helped to analyze the mathematical structure of the problem and to derive its decomposition in the form of a 1-rank tensor decomposition known as Canonical Decomposition (CANDECOMP) \cite{Harshman1970,KoldaTensorReview2009,AcarTensorSurvey2009,Martin2004}. More details about proposed tensor formulation and solving algorithm can be found in our companion paper.

The identified tensorial structure allowed efficient computation of the tensor multiplication which therefore can be efficiently used within a tensor Krylov solver. Another important implication of the identified property of the tensor formulation, is a dramatic reduction of storage requirements. Our approach does not require explicit storage either of the system tensor or its decomposition, and thus the problem of size, discussed above (33'554'432 unknowns parameters), for about 9'000'000 measurements, can be computed on current inexpensive multi core computer equipped with only 2 GBytes of physical memory in less than 60 minutes, surpassing the capability of a published, matrix-based solving algorithm \cite{ArigovindanImageRec2005} for the same problem by far. The decomposability of the tensor multiplication allows convenient and \textbf{\textit{complete parallelization}} of computations using all currently available parallelism paradigms, including Single Instruction Multiple Data (SIMD), Multi Core, Clusters of PCs and General Purpose GPU (GPGPU) computing technology.
	
	Note that parallelizability of tensor operations, in general, has proven to be very beneficial: on Multi Core CPU's we have observed a speed increase almost proportional to the number of cores; and on a recent graphics card (ATI Radeon HD 4870 X2) an additional increase has been measured for specific examples, of more than 10, compared to a Quad Core CPU.

\section{Discussion}
\label{sec5}
In this work we introduced a framework for computational solving of large tensor structured problems. The proposed approach leads to a natural, dimensionality-preserving formulation of tensor structured problems, and by maintaining the multidimensional structure of the data, and the commutativity properties of the framework, directly induces computationally efficient solving strategies that can profit from automatic expression analysis, separability properties of the tensor formulation, and the preservation of spatial coherence of the data, to speed up convergence in tensor solvers as compared to their matrix counterparts.		
	
The chosen formalism based on Einstein notation together with introduced commutativity of tensor multiplication, makes the framework well suited for implementing a tensor programming language, offering self-optimized computations of tensor expressions by semantic analysis of their terms while avoiding explicit use of tensor indexing in the actual implementations of solving programs. This was verified by implementing our own tensor library built on the principles of the proposed framework. 

The properties presented above distinguish our framework from one often used in the field of tensor decompositions \cite{Tucker1966,Leibovici1998,LathauwerHOSVD2000,Harshman1970,KoldaTensorReview2009,AcarTensorSurvey2009,Martin2004}, but do not limit its use for solving problems studied in that field.

The ongoing research shows that the framework is well suited for solving problems of signal processing in higher dimensions, which have inherent tensor structure. We are currently investigating its usefulness in problems from image reconstruction, computer vision and fluid dynamics.

In summary, this high-level approach fits well into the landscape of increasing interest in tensors as a workhorse for elegant and efficient solving of those scientific problems which arise from multidimensional data.

\footnotetext[4]{Objects C,D,E,F,G,H correspond respectively to $\cal{C},\cal{D},\cal{E},\cal{F},\cal{G},\cal{H}$ from tensor expressions; "*" denotes the tensor product}

\footnotetext[5]{The tensor representation of finite difference is the base for derivative related operators such as gradient, divergence, curl and Laplacian.}
	
\begin{table}[ht]
\centering \footnotesize
\caption{Basic multidimensional data processing building blocks in separable tensor form}
\begin{tabular}{|l|l|l|} \hline

Signal processing operator & Tensor expression & Corresponding \\ & & implementation \footnotemark[4] \\ \hline

Discrete convolution & ${\cal{F}}^{x_1 y_1 z_1} = {\cal{E}}^{x_1}_{x} \cdot {\cal{G}}^{y_1}_{y} \cdot {\cal{H}}^{z_1}_{z} \cdot {\cal{C}}^{xyz}$ & F := E*G*H*C; \\
& ${\cal{E,G,H}}$ - 1D convolution transforms & \\ \hline

Finite difference \footnotemark[5] & ${\cal{F}}^{x_1 y z} = {\cal{E}}^{x_1}_{x} \cdot {\cal{C}}^{xyz}$ & F := E*C; \\
& ${\cal{F}}^{x y_1 z} = {\cal{G}}^{y_1}_{y} \cdot {\cal{C}}^{xyz}$ & F := G*C; \\
& ${\cal{F}}^{x y z_1} = {\cal{H}}^{z_1}_{z} \cdot {\cal{C}}^{xyz}$ & F := H*C; \\
& ${\cal{E,G,H}}$ - 1D finite difference transforms &	\\ \hline

Discrete Fourier transform & ${\cal{F}}^{x_1 y_1 z_1} = {\cal{E}}^{x_1}_{x} \cdot {\cal{G}}^{y_1}_{y} \cdot {\cal{H}}^{z_1}_{z} \cdot {\cal{C}}^{xyz}$ & F := E*G*H*C; \\
& ${\cal{E,G,H}}$ - 1D DFT transforms &	\\ \hline

Rotation & ${\cal{F}}^{x_2 y_1 z_1} = {\cal{D}}^{x_2}_{x_1} \cdot \left( {\cal{E}}^{x_1}_{x} \cdot {\cal{G}}^{y_1}_{y} \cdot {\cal{H}}^{z_1}_{z} \right) \cdot {\cal{C}}^{xyz}$ & F := D*E*G*H*C; \\
& ${\cal{D,E,G,H}}$ - 1D shear transforms & \\ \hline

Upsampling/downsampling & ${\cal{F}}^{u v w} = {\cal{E}}^{u}_{x} \cdot {\cal{G}}^{v}_{y} \cdot {\cal{H}}^{w}_{z} \cdot {\cal{C}}^{xyz}$ & F := E*G*H*C; \\
& ${\cal{E,G,H}}$ - 1D up-/downsampling transforms & \\ \hline
\end{tabular}
\end{table}

\begin{table}[ht]
\centering \footnotesize
\caption{Large scientific problems suited for solution by computational Tensor Algebra}
\begin{tabular}{|l|l|} \hline
Computational Application & Scientific Field \\ \hline
Multilinear data models using tensor decompositions & Signal Processing \\
(e.g. CANDECOMP, TUCKER, HOSVD): & Machine Learning \\

- Pattern recognition \cite{SavasDigitClassHOSVD2007,VasilescuTensorFaces2002}	& Computer Vision \\
- Data compression \cite{WangDataComprCANDECOMP2004} & Data Mining \\
- Independent Component Analysis \cite{ComonTensorDecomp2000,Lathauwer1998} & Chemistry \\
- Signal filtering \cite{Muti2003} & Neuroscience \\
- Modeling of fluorescence data \cite{AcarTensorSurvey2009} & \\
- Social Network Analysis \cite{AcarTensorSurvey2009} & \\ \hline

Solving inverse multidimensional problems: & Signal Processing \\
- Least squares problems \cite{ArigovindanImageRec2005,Pereyra1973} & Computer Vision \\
- Differential equations (Poisson, Navier-Stokes, Finite & Machine Learning \\
Element Methods) \cite{Lynch1964,GavrilyukHierarchTensor2005} & Fluid Dynamics \\
- Integral equations \cite{EldenTensorLS2005} & Chemistry \\
& Electromagnetism \\ \hline

Modeling of properties of complex systems consisting of many  & Computational Physics \\
interacting elements expressed by tensor expressions: & Quantum Chemistry \\
- Modeling of electronic and optical properties of & Computational Chemistry\\
molecules and their interactions \cite{Schmidt1993,Raghavachari1991} & Material Science \\ \hline

\end{tabular}
\end{table}

\appendix
\section{APPENDIX}

\subsection{Dual vectors spaces, contravariant and covariant mechanisms of transformation}
\label{cov_contrav}
Let $V$ be a real vector space with finite number of dimensions $N_V$. We designate by $V^*$ unique vector space of same dimensionality that is dual in respect to $V$. Any element $\mathbf{t}^* \in V^*$ is a linear map $V \to \Re$

\begin{equation}
\langle \mathbf{t}^* , \mathbf{t}  \rangle \equiv \mathbf{t}^* (\mathbf{t}) \in \Re
\label{A1}
\end{equation}

For a given base $e_v$ of $V$ transformations to a new base $\hat{e}_{v_1}$ and vice versa are defined by

\begin{equation}
\hat{e}_{v_1} = \sum_{v=1}^{N_V}{{\cal{P}}^{v}_{v_1} \cdot e_v}, \: e_{v_2} = \sum_{v_1=1}^{N_V}{{\cal{S}}^{v_1}_{v_2} \cdot \hat{e}_{v_1}}
\label{A2}
\end{equation}

where ${\cal{P}}^{v}_{v_1}$ and ${\cal{S}}^{v_1}_{v_2}$ are direct and indirect base transformations respectively. A vector from $V$ given by $\sum_{v=1}^{N_V}{{\cal{T}}^{v} \cdot e_v}$ according to (\ref{A2}) in a new coordinate system is represented by components

\begin{equation}
\hat{{\cal{T}}}^{v_1} = \sum_{v=1}^{N_V}{\cal{S}}^{v_1}_{v} \cdot {\cal{T}}^{v}
\label{A3}
\end{equation}

Components of such vectors which transform indirectly in respect to transformation of the base are called contravariant.

For the base $e_v$ of $V$ and dual base $e^{v_1}$ of $V^*$ the following equation is satisfied

\begin{equation}
\langle e^{v_1} , e_v \rangle \equiv \delta^{v_1}_{v} = 
		\left\{ 
			\begin{array}{ll}
         		1, & \mbox{if $v_1 = v$}\\
         		0, & \mbox{if $v_1 \neq v$} 
         	\end{array} 
         \right.
\label{A4}         
\end{equation}

where $\delta^{v_1}_{v}$ is Kronecker delta. Thus the dual base $e^{v_1}$ transforms according to

\begin{equation}
\hat{e}^{v_2} = \sum_{v_1=1}^{N_V}{{\cal{S}}^{v_2}_{v_1}} \cdot e^{v_1}, \: 
e^{v_3} = \sum_{v_2=1}^{N_V}{{\cal{P}}^{v_3}_{v_2} \cdot \hat{e}^{v2}}
\label{A5}
\end{equation}

In contrast to contravariant vector components, components of vectors from dual space $V^*$ transform in the same way as base of $V$

\begin{equation}
{\cal{T}}_{v_2} = \sum_{v_1=1}^{N_V}{{\cal{P}}^{v_1}_{v_2} \cdot {\cal{T}}_{v_1}}
\label{A6}
\end{equation}

Components of such vectors are called covariant.

\subsection{Generalized definition of a tensor}

\definition Let ordered set of real vector spaces $U_1, \ldots , U_P$ , $V_1, \ldots , V_Q$ with finite dimensions $I_1, \ldots , I_P$ , $J_1, \ldots , J_Q$ and their dual vector space $U^*_1, \ldots , U^*_P$ , $V^*_1, \ldots , V^*_{Q}$ be given. By identifying $(U^*_i)^*$ with $U_i$ we define every vector $\mathbf{u}_i \in U_i$ to be a linear map $U^*_i \to \Re$ given by

\begin{equation}
\mathbf{u}_i(\mathbf{u}^*_i) \equiv \langle \mathbf{u}^*_i , \mathbf{u}_i \rangle \in \Re
\label{A7}
\end{equation}

for arbitrary $\mathbf{u}^*_i \in U^*_i$, where $\langle , \rangle$ denotes the inner product. With $(P+Q)$ vectors $\mathbf{u}_1 \in U_1, \ldots , \mathbf{u}_P \in U_P$, $\mathbf{v}^*_1 \in V^*_1, \ldots , \mathbf{v}^*_Q \in V^*_Q$ let the element denoted ${\cal{T}} = \mathbf{u}_1 \times \cdots \times \mathbf{u}_P \times \mathbf{v}^*_1 \times \cdots \times \mathbf{v}^*_Q$ be a  $(P+Q)$-linear map from $U_1 \times \cdots \times U_P \times V^*_1 \times \cdots \times V^*_Q$ to $\Re$ defined by

\begin{equation}
\begin{split}
\mathbf{u}_1 \times \cdots \times \mathbf{u}_P \times \mathbf{v}^*_1 \times \cdots \times \mathbf{v}^*_Q \left(
\mathbf{a}^*_1, \ldots , \mathbf{a}^*_P, \mathbf{b}_1, \ldots , \mathbf{b}_Q
\right) =\\
\langle \mathbf{a}^*_1 , \mathbf{u}_1 \rangle \cdots \langle \mathbf{a}^*_P , \mathbf{u}_P \rangle 
\langle \mathbf{v}^*_1 , \mathbf{b}_1 \rangle \cdots \langle \mathbf{v}^*_Q , \mathbf{b}_Q \rangle
\end{split}
\label{A8}
\end{equation}

where $\mathbf{a}^*_i$, $\mathbf{b}_j$ are arbitrary vectors in $U^*_i$ and $V_j$ respectively. The element $\cal{T}$ is termed a  $(P+Q)$-th order decomposed tensor,  $P$-times contravariant and  $Q$-times covariant. The space generated by all linear combinations of decomposed tensors is termed the tensor product space. Any arbitrary tensor can be represented as a weighted sum of decomposed tensors \cite{Harshman1970,AcarTensorSurvey2009,Martin2004,Lathauwer1998}. This definition is compatible with, and extends \cite{Leibovici1998,Pereyra1973}.

\subsection{Tensor notation}
\label{Tensor_notation}
To each dimension of the physical space we assign a vector space with a given dimensionality. For example for 3-dimensional physical space we define vector spaces $X,Y,Z$ with sizes $N_X, N_Y, N_Z$. In such vector spaces one can define multiple coordinate systems related to each other by transformations discussed above. Components of a tensor in some tensor product space are designated as an array with multiple indices. Each index in the designation has direct correspondence to a dimension of the physical space and thus to one of the defined vector spaces. For example for 3-dimensional case a tensor which is an element of $X \times Y \times Z$ can be designated by ${\cal{T}}^{x_i y_j z_k}$. Subindices $i$, $j$, $k$ are used for distinguishing between different coordinate systems.

\subsection{Elementwise tensor operations}
\label{tensor_operations}
Linearity of the defined tensor product space implies definition of the following elementwise tensor operations such as tensor addition, subtraction and multiplication by a scalar:

\begin{align}
&{\cal{A}}^{x_1 y_1 z_1}_{xyz} + {\cal{B}}^{x_1 y_1 z_1}_{xyz} = {\cal{C}}^{x_1 y_1 z_1}_{xyz} \\
&{\cal{A}}^{x_1 y_1 z_1}_{xyz} - {\cal{B}}^{x_1 y_1 z_1}_{xyz} = {\cal{D}}^{x_1 y_1 z_1}_{xyz} \\
&\lambda \cdot {\cal{A}}^{x_1 y_1 z_1}_{xyz} = {\cal{E}}^{x_1 y_1 z_1}_{xyz}
\label{A9}
\end{align}

With introduced constraint of ordered vector spaces all properties of these operations such as associativity and commutativity remain unchanged.

\subsection{Kronecker delta}
\label{kronecker_delta}
Kronecker delta presented in expression (\ref{A4}) is a second order tensor which is equivalent to the identity transformation. Mixed version of the tensor applied to components of a contravariant or covariant first order tensor (vector) does not introduce changes

\begin{align}
&\delta^{x_1}_{x} \cdot {\cal{T}}^x = {\cal{T}}^{x_1} \equiv {\cal{T}}^{x} \\
&\delta^{x_1}_{x} \cdot {\cal{U}}_{x_1} = {\cal{U}}_{x} \equiv {\cal{U}}_{x_1}
\label{A10}
\end{align}

Covariant or contravariant versions of the Kronecker delta change variance of tensor indices, but do not change values of tensor components

\begin{align}
&\delta_{{x_1}{x}} \cdot {\cal{T}}^x = {\cal{T}}_{x_1} \equiv {\cal{T}}_{x} \\
&\delta^{{x_1}{x}} \cdot {\cal{U}}_{x_1} = {\cal{U}}^{x} \equiv {\cal{U}}^{x_1}
\label{A11}
\end{align}

Note that here indices $x$ and $x_1$ correspond to identical coordinate systems. In the same way Kronecker delta can be applied along some dimension of a higher order tensor

\begin{equation}
\delta^{x_1}_{x} \cdot {\cal{T}}^{xyz} \equiv \delta^{y_1}_{y} \cdot {\cal{T}}^{xyz} \cdot \delta^{z_1}_{z} \cdot {\cal{T}}^{xyz} \equiv {\cal{T}}^{xyz}
\label{A11}
\end{equation}

The tensor multiplication of Kronecker deltas for different vector spaces forms multidimensional identity transformation

\begin{align}
&\delta^{x_1 y_1 z_1}_{xyz} = \delta^{x_1}_{x} \cdot \delta^{y_1}_{y} \cdot \delta^{z_1}_{z} \\
&\delta^{x_1 y_1 z_1}_{xyz} \cdot {\cal{T}}^{xyz} = {\cal{T}}^{x_1 y_1 z_1} \equiv {\cal{T}}^{xyz}
\label{A12}
\end{align}

\subsection{Inner product}
\label{InnerProd}
the inner product of two third order Euclidean tensors ${\cal{A}}^{xyz}$ and ${\cal{B}}^{xyz}$ is defined by

\begin{align}
&\langle A, B \rangle = \left( {\cal{A}}^{xyz} \cdot {\delta}_{x{x_1}y{y_1}z{z_1}} \right) \cdot \left( {\cal{B}}^{xyz} \cdot {\delta}^{{x_1}{y_1}{z_1}}_{xyz} \right) = \hat{{\cal{A}}}^{{x_1}{y_1}{z_1}} \cdot \hat{{\cal{B}}}_{{x_1}{y_1}{z_1}}
\label{A13}
\end{align}

At the same time using the language of elementwise products the inner product can be expressed as

\begin{align}
&\langle A, B \rangle = \left( {\cal{A}}^{xyz} \cdot {\cal{B}}^{xyz} \right) \cdot \left({\delta}_{x{x_1}y{y_1}z{z_1}}  \cdot {\delta}^{{x_1}{y_1}{z_1}}_{xyz} \right) = {\cal{C}}^{{x_1}{y_1}{z_1}} \cdot {\cal{I}}_{{x_1}{y_1}{z_1}}
\label{A14}
\end{align}

where ${\cal{I}}_{{x_1}{y_1}{z_1}}$ is a tensor with all components equal to one, which can be easily checked based on the properties of Kronecker delta.

Note that in general case of non-Euclidean spaces the inner product is defined with the use of the metric tensor which for Euclidean case reduces to the Kronecker delta.

\subsection{Tensor transposition}
\label{tensor_transp}
In classical settings transposition of a tensor is defined by changing positions of tensor indices. But in our framework the order of indices is unique according to a predefined order of vector spaces. Thus patterns like ${\textbf{A}^T}\textbf{A}$ from matrix normal equations correspond to the tensor multiplication with change of variance by the Kronecker delta in Euclidean spaces or the metric tensor in general case.

\subsection{Differentiation in respect to a tensor}
\label{tensor_diff}
When dealing with an optimization problem, differentiation in respect to the optimization parameter may be required. Consider an example of differentiating a tensor valued function of the form  
$f({\cal{U}}) = {\cal{A}}^{x_1 y_1 z_1}_{xyz} \cdot {\cal{U}}^{xyz}$

\begin{equation}
{\partial \over {\partial {\cal{U}}}} f({\cal{U}}) = {\partial \over {\partial {\cal{U}}^{x_2 y_2 z_2}}} \left(
{\cal{A}}^{x_1 y_1 z_1}_{xyz} \cdot {\cal{U}}^{xyz}
\right) = {\cal{A}}^{x_1 y_1 z_1}_{xyz} \cdot \left( {\partial {\cal{U}}^{xyz}} \over {\partial {\cal{U}}^{x_2 y_2 z_2}} \right)
\label{A15}
\end{equation}

where indices $(x, y, z)$ and $(x_2, y_2, z_2)$ correspond to identical coordinate systems. Using tensor transformational properties it is straightforward to show that

\begin{equation}
{{\partial {\cal{U}}^{xyz}} \over {\partial {\cal{U}}^{x_2 y_2 z_2}}} = \delta^{xyz}_{x_2 y_2 z_2}
\label{A16}
\end{equation}

Thus we have

\begin{equation}
{\partial \over {\partial {\cal{U}}}} f({\cal{U}}) = {\cal{A}}^{x_1 y_1 z_1}_{xyz} \cdot \delta^{xyz}_{x_2 y_2 z_2} = {\cal{A}}^{x_1 y_1 z_1}_{x_2 y_2 z_2} \equiv {\cal{A}}^{x_1 y_1 z_1}_{xyz}
\label{A17}
\end{equation}

In a similar way one can compute derivative of more complex scalar and tensor valued functions.

%
\bibliographystyle{spmpsci}       
\bibliography{references}         

\end{document}